\documentclass[useAMS]{article}

\usepackage{graphicx}
\begin{document}
\title{{\itshape Nearly Pentadiagonal Linear Systems}}

\author{ A. A.KARAWIA\footnote{ Home Address: Mansoura University,
Mansoura 35516, Egypt. E-mail: abibka@mans.edu.eg} $\quad$and S.A.
El-Shehawy\footnote{ Home Address: Menoufia University, Shebin
El-Kom, Egypt. E-mail addresses: shshehawy64@yahoo.com}\\
\\
Qassim University, PYP, Buraidah, 51452, Saudi Arabia\\
  Qassim University, College of Science, Department of Mathematics, P.O. Box 237, \\Buriedah 51452, Saudi Arabia}

\date{}

\makeatletter \@addtoreset{equation}{section}
\renewcommand{\theequation}{\thesection.\@arabic\c@equation}
\makeatother

\maketitle
\begin{abstract}
In this paper we present  efficient computational and symbolic
algorithms for solving  a nearly pentadiagonal linear systems. The
implementation of the algorithms using Computer Algebra Systems
(CAS) such as MAPLE, MACSYMA, MATHEMATICA, and MATLAB is
straightforward. Two examples are given in order to illustrate the
algorithms. \bigskip

\begin{flushleft}\footnotesize
{\textbf{Keywords:} pentadiagonal matrices; nearly pentadiagonal
matrices; linear systems; determinants; computer algebra systems
(CAS).\bigskip

}
\end{flushleft}

\newtheorem{alg}{Algorithm}[section]

\end{abstract}

\section{Introduction}

Many problems in mathematics and applied science require the
solution of linear systems having nearly pentadiagonal coefficient
matrices. This kind of linear system arises in many fields of
numerical computation and differential Equations [1, 2,3]. This article is a general case of the author article [3].\\

The main goal of the current paper is to develop an efficient
algorithms for solving a general nearly pentadiagonal linear systems
of the form:
\begin{equation}
AX=Y
\end{equation}
where
\begin{eqnarray}
A=\left[\begin{array}{ccccccccc}
                     d_1 & a_1 & \tilde{a}_1 & s & 0 & 0 & 0 & \ldots & 0 \\
                     b_2 & d_2 & a_2 & \tilde{a}_2 & 0 & 0 & 0 & \ldots & 0 \\
                     \tilde{b}_3 & b_3 & d_3 & a_3 & \tilde{a}_3 & 0 & 0 & \ldots & \vdots \\
                     0 & \tilde{b}_4 & b_4 & d_4 & a_4 & \tilde{a}_4 & 0 & \ldots & \vdots \\
                     \vdots & \ddots & \ddots & \ddots & \ddots & \ddots & \ddots & \ddots & \vdots\\
                     \vdots & \ddots & \ddots & \ddots & \ddots & \ddots & \ddots & \ddots & 0 \\
                     \vdots & \ldots & \ldots & 0 & \tilde{b}_{n-2} & b_{n-2} & d_{n-2} & a_{n-2} & \tilde{a}_{n-2} \\
                     \vdots & \ldots & \ldots & \ldots & 0 & \tilde{b}_{n-1} & b_{n-1} & d_{n-1} & a_{n-1} \\
                     0 & \ldots & \ldots & \ldots  & 0 & t & \tilde{b}_n & b_n & d_n
                   \end{array}
        \right],
\end{eqnarray}
$X=(x_1,x_2,\ldots,x_n)^T$, $Y=(y_1,y_2,\ldots,y_n)^T$ and $n\ge5$.\\
\\
A general $n\times n$ nearly pentadiagonal matrix $A$ of the form
(1.2) can be stored in $5n-4$ memory locations by using five vectors
$\mbox{\boldmath$\tilde{a}$}=(\tilde{a}_1,\tilde{a}_2,\ldots,\tilde{a}_{n-2})$,
$\mbox{\boldmath$a$}=(a_1,a_2,\ldots,a_{n-1},s)$,
$\mbox{\boldmath$b$}=(t,b_2,b_3,\ldots,b_n)$,
$\mbox{\boldmath$\tilde{b}$}=(\tilde{b}_3,\tilde{b}_4,\ldots,\tilde{b}_n)$,
and $\mbox{\boldmath$d$}=(d_1,d_2,\ldots,d_n)$. When considering the
system (1.1) it is advantageous to introduce three additional
vectors $\mbox{\boldmath$c$}=(c_1,c_2,\ldots,c_n)$,
$\mbox{\boldmath$e$}=(e_1,e_2,\ldots,e_n)$ and
$\mbox{\boldmath$f$}=(f_1,f_2,f_3,\ldots,f_n)$. These vectors are
related to the vectors $\mbox{\boldmath$\tilde{a}$}$,
$\mbox{\boldmath$a$}$, $\mbox{\boldmath$d$}$, $\mbox{\boldmath$b$}$, and $\mbox{\boldmath$\tilde{b}$}$.\\

 The current paper is organized as follows. In section 2, the main
 results are given. Illustrative examples are presented in section 3. In section 4, a
 conclusion is given.

\section{Main results}
In this section we are going to formulate a new computational and
symbolic algorithms for solving a general nearly pentadiagonal
linear systems of the form (1.1). To do this, we begin by
considering the $LU$ decomposition [4] of the matrix $A$
in the form:\\
\begin{eqnarray}
\left[\begin{array}{ccccccccc}
                     d_1 & a_1 & \tilde{a}_1 & s & 0 & 0 & 0 & \ldots & 0 \\
                     b_2 & d_2 & a_2 & \tilde{a}_2 & 0 & 0 & 0 & \ldots & 0 \\
                     \tilde{b}_3 & b_3 & d_3 & a_3 & \tilde{a}_3 & 0 & 0 & \ldots & \vdots \\
                     0 & \tilde{b}_4 & b_4 & d_4 & a_4 & \tilde{a}_4 & 0 & \ldots & \vdots \\
                     \vdots & \ddots & \ddots & \ddots & \ddots & \ddots & \ddots & \ddots & \vdots\\
                     \vdots & \ddots & \ddots & \ddots & \ddots & \ddots & \ddots & \ddots & 0 \\
                     \vdots & \ldots & \ldots & 0 & \tilde{b}_{n-2} & b_{n-2} & d_{n-2} & a_{n-2} & \tilde{a}_{n-2} \\
                     \vdots & \ldots & \ldots & \ldots & 0 & \tilde{b}_{n-1} & b_{n-1} & d_{n-1} & a_{n-1} \\
                     0 & \ldots & \ldots & \ldots  & 0 & t & \tilde{b}_n & b_n & d_n
                   \end{array}
        \right]=\nonumber
\end{eqnarray}
\begin{eqnarray}
          \left[\begin{array}{ccccccccc}
                     1 & 0 & 0 & 0 & 0 & 0 & 0 & \ldots & 0 \\
                     f_2 & 1 & 0 & 0 & 0 & 0 & 0 & \ldots & 0 \\
                     \frac{\tilde{b}_3}{c_1} & f_3 & 1 & 0 & 0 & 0 & 0 & \ldots & \vdots \\
                     0 & \frac{\tilde{b}_4}{c_2} & f_4 & 1 & 0 & 0 & 0 & \ldots & \vdots \\
                     \vdots & \ddots & \ddots & \ddots & \ddots & \ddots & \ddots & \ddots & \vdots\\
                     \vdots & \ddots & \ddots & \ddots & \ddots & \ddots & \ddots & \ddots & 0 \\
                     \vdots & \ldots & \ldots & 0 & \frac{\tilde{b}_{n-2}}{c_{n-4}} & f_{n-2} & 1 & 0 & 0 \\
                     \vdots & \ldots & \ldots & \ldots & 0 & \frac{\tilde{b}_{n-1}}{c_{n-3}} & f_{n-1} & 1 & 0 \\
                     0 & \ldots & \ldots & \ldots  & 0 & \frac{t}{c_{n-3}} & f_1 & f_n &
                     1
                   \end{array}
        \right]\left[\begin{array}{ccccccccc}
                     c_1 & e_1 & \tilde{a}_1 & s & 0 & 0 & 0 & \ldots & 0 \\
                     0 & c_2 & e_2 & e_n & 0 & 0 & 0 & \ldots & 0 \\
                     0 & 0 & c_3 & e_3 & \tilde{a}_3 & 0 & 0 & \ldots & \vdots \\
                     0 & 0 & 0 & c_4 & e_4 & \tilde{a}_4 & 0 & \ldots & \vdots \\
                     \vdots & \ddots & \ddots & \ddots & \ddots & \ddots & \ddots & \ddots & \vdots\\
                     \vdots & \ddots & \ddots & \ddots & \ddots & \ddots & \ddots & \ddots & 0 \\
                     \vdots & \ldots & \ldots & 0 & 0 & 0 & c_{n-2} & e_{n-2} & \tilde{a}_{n-2} \\
                     \vdots & \ldots & \ldots & \ldots & 0 & 0 & 0 & c_{n-1} & e_{n-1} \\
                     0 & \ldots & \ldots & \ldots  & 0 & 0 & 0 & 0 & c_n
                   \end{array}
        \right]
\end{eqnarray}
From (2.1) we obtain
\begin{equation}
c_i= \left\{\begin{array}{ll}
d_1 & \textrm{if $i=1$}\\
d_2-f_2 e_1 & \textrm{if $i=2$}\\
d_i-f_i
e_{i-1}-\frac{\tilde{b}_i}{c_{i-2}}\tilde{a}_{i-2},\tilde{a}_{2}=e_n
& \textrm{if $i=3(1)n-1$}\\
d_n-f_1\tilde{a}_{n-2}-f_n e_{n-1}& \textrm{if $i=n$},
\end{array}\right.
\end{equation}
where
\\
\begin{equation}
e_i= \left\{\begin{array}{ll}
a_1 & \textrm{if $i=1$}\\
a_2-f_2\tilde{a}_1 & \textrm{if $i=2$}\\
\tilde{a}_2-f_2 a_n & \textrm{if $i=n$}\\
a_3-\frac{\tilde{b}_3}{c_1}a_n-f_3 e_n & \textrm{if $i=3$}\\
a_i-f_i\tilde{a}_{i-1} & \textrm{if $i=4(1)n-1$},
\end{array}\right.
\end{equation}
\\
and
\begin{equation}
\hspace*{-2.7cm}f_i= \left\{\begin{array}{ll}
(\tilde{b}_n-\frac{t}{c_{n-3}}e_{n-3})/c_{n-2} & \textrm{if $i=1$}\\
\frac{b_2}{c_1} & \textrm{if $i=2$}\\
(b_i-\frac{\tilde{b}_i}{c_{i-2}})/c_{i-1} & \textrm{if $i=3(1)n-1$}\\
(b_n-\frac{t}{c_{n-3}}\tilde{a}_{n-3}-f_1 e_{n-2})/c_{n-1} &
\textrm{if $i=n$ }.
\end{array}\right.
\end{equation}
\\
It is not difficult to prove that the $LU$ decomposition (2.1)
exists only if $c_i\ne 0,\quad i=1(1)n-1$. Moreover, a general nearly
pentadiagonal linear system (1.1) possesses a unique solution if, in
addition, $c_n\ne 0$. On the other hand, the determinant of the
matrix $A$ is given by:
\begin{equation}
det(A)=\prod_{i=1}^n c_i,
\end{equation}
and this shows the importance of the vector
$\mbox{\boldmath$c$}$ [5].\\
\\
We may now formulate the following results.
\begin{alg}
To solve the general nearly pentadiagonal linear system (1.1),
we may proceed as follows:\\
\textbf{step 1:} Set $c_1=d_1$.\\
\textbf{step 2:} If  $c_1=0$, then OUTPUT('the method is fails'); STOP.\\
\textbf{step 3:} Set $a_n=s$, $b_1=t$, $e_1=a_1,$
$f_2=\frac{b_2}{c_1}$, and $c_2=d_2-f_2 e_1$.\\
\textbf{step 4:} If  $c_2=0$, then OUTPUT('the method is fails'); STOP.\\
 \textbf{step 5:} Set $e_2=a_2-f_2\tilde{a}_1$,
 $e_n=\tilde{a}_2-f_2a_n$, and $\tilde{a}_2=e_n$.\\
\textbf{step 6:} For $i=3, 4,\ldots,n-1$ Compute\\
\hspace*{2.3cm}$f_i=(b_i-\frac{\tilde{b}_i}{c_{i-2}}e_{i-2})/c_{i-1}$,\\
\hspace*{2.3cm}$c_i=d_i-f_ie_{i-1}-\frac{\tilde{b}_i}{c_{i-2}}\tilde{a}_{i-2}$,\\
\hspace*{2.3cm}If  $c_i=0$, then OUTPUT('the method is fails'); STOP,\\
\hspace*{2.3cm}If $i=3$, then $e_i=a_i-f_i\tilde{a}_{i-1}-\frac{\tilde{b}_i}{c_{i-2}}a_n$ else $e_i=a_i-f_i\tilde{a}_{i-1}$.\\
\textbf{step 7:} Compute\\
\hspace*{2.3cm}$f_1=(\tilde{b}_n-\frac{b_1}{c_{n-3}}e_{n-3})/c_{n-2}$,\\
\hspace*{2.3cm}$f_n=(b_n-\frac{t}{c_{n-3}}\tilde{a}_{n-3}-f_1e_{n-2})/c_{n-1}$,\\
\hspace*{2.3cm}$c_n=d_n-f_1\tilde{a}_{n-2}-f_ne_{n-1}$,\\
\textbf{step 8:} Set $z_1=y_1$, $z_2=y_2-f_2 z_1$.\\
\textbf{step 9:} For $i=3,4, \ldots,n-1$ Compute\\
\hspace*{2.3cm}$z_i=y_i-f_iz_{i-1}-\frac{\tilde{b}_i}{c_{i-2}}z_{i-2}$.\\
\textbf{step 10:} Set
$z_n=y_n-f_nz_{n-1}-f_1z_{n-2}-\frac{t}{c_{n-3}}z_{n-3}$.\\
\textbf{step 11:} Compute the solution vector
$\mbox{\boldmath$x$}$ using\\
\hspace*{2cm} $x_n=\frac{z_n}{c_n}$,
$x_{n-1}=\frac{z_{n-1}-e_{n-1} x_n}{c_{n-1}}$,\\
\hspace*{2cm}For $i=n-2,n-3,\ldots, 2$
compute\\
\hspace*{2.9cm}$x_i=\frac{z_i-e_i x_{i+1}-\tilde{a}_i x_{i+2}}{c_i}$,\\
\hspace*{2cm}Set $x_1=\frac{z_1-e_1 x_2-\tilde{a}_1 x_3-s x_4}{c_1}$.\\

\end{alg}
The new algorithm 2.1 will be referred to as \textbf{KNPENTA}
algorithm. \textbf{KNPENTA} algorithm for solving the nearly
pentadiagonal system (1.1) is generally preferable because the
conditions $c_i\ne 0,\quad i=1(1)n$ are sufficient for its validity. The advantage of the vector $\mbox{\boldmath$c$}$ is now
clear.\\

The following symbolic algorithm is developed in order to remove the
cases where the numeric algorithm \textbf{KNPENTA} fails.

\begin{alg}
To solve the general nearly pentadiagonal linear system (1.1),
we may proceed as follows:\\
\textbf{step 1:} Set $c_1=d_1$.\\
\textbf{step 2:} Set  $c_1=x$ and $d_1=x$($x$ is just a symbolic name) whenever $c_1=0$.\\
\textbf{step 3:} Set $a_n=s$, $b_1=t$, $e_1=a_1$,
$f_2=\frac{b_2}{c_1}$, and $c_2=d_2-f_2 e_1$.\\
\textbf{step 4:} Set  $c_2=x$ whenever $c_2=0$.\\
\textbf{step 5:} Set $e_2=a_2-f_2\tilde{a}_1$,
$e_n=\tilde{a}_2-f_2a_n$, and $\tilde{a}_2=e_n$.\\
\textbf{step 6:} For $i=3, 4,\ldots,n-1$ Compute\\
\hspace*{2.3cm}$f_i=(b_i-\frac{\tilde{b}_i}{c_{i-2}}e_{i-2})/c_{i-1}$,\\
\hspace*{2.3cm}$c_i=d_i-f_ie_{i-1}-\frac{\tilde{b}_i}{c_{i-2}}\tilde{a}_{i-2}$,\\
\hspace*{2.3cm}Set  $c_i=x$ whenever $c_i=0$,\\
\hspace*{2.3cm}If $i=3$, then $e_i=a_i-f_i\tilde{a}_{i-1}-\frac{\tilde{b}_i}{c_{i-2}}a_n$ else $e_i=a_i-f_i\tilde{a}_{i-1}$.\\
\textbf{step 7:} Compute\\
\hspace*{2.3cm}$f_1=(\tilde{b}_n-\frac{b_1}{c_{n-3}}e_{n-3})/c_{n-2}$,\\
\hspace*{2.3cm}$f_n=(b_n-\frac{t}{c_{n-3}}\tilde{a}_{n-3}-f_1e_{n-2})/c_{n-1}$,\\
\hspace*{2.3cm}$c_n=d_n-f_1\tilde{a}_{n-2}-f_ne_{n-1}$,\\
\hspace*{2.3cm}Set  $c_n=x$ whenever $c_n=0$.\\
\textbf{step 8:} Set $z_1=y_1$, $z_2=y_2-f_2 z_1$.\\
\textbf{step 9:} For $i=3,4, \ldots,n-1$ Compute\\
\hspace*{2.3cm}$z_i=y_i-f_iz_{i-1}-\frac{\tilde{b}_i}{c_{i-2}}z_{i-2}$.\\
\textbf{step 10:} Set
$z_n=y_n-f_nz_{n-1}-f_1z_{n-2}-\frac{t}{c_{n-3}}z_{n-3}$.\\
\textbf{step 11:} Compute the solution vector
$\mbox{\boldmath$x$}$ using\\
\hspace*{2cm} $x_n=\frac{z_n}{c_n}$,
$x_{n-1}=\frac{z_{n-1}-e_{n-1} x_n}{c_{n-1}}$,\\
\hspace*{2cm}For $i=n-2,n-3,\ldots, 2$
compute\\
\hspace*{2.9cm}$x_i=\frac{z_i-e_i x_{i+1}-\tilde{a}_i x_{i+2}}{c_i}$,\\
\hspace*{2cm}Set $x_1=\frac{z_1-e_1 x_2-\tilde{a}_1 x_3-s x_4}{c_1}$.\\
\textbf{step 12:} Substitute $x=0$ in all expressions of the
solution vector $x_i,i=1,2,\ldots,n$.\\
\end{alg}

The symbolic algorithm 2.2 will be referred to as \textbf{KSNPENTA}
algorithm.\\
In [6], Claerbout showed that the two-dimensional Laplacian
operator, which appears in 3-D finite-difference migration, has the
form of pentadiagonal matrix. If we choose $d_i=-4,\quad
i=1(1)n$,$\quad s=t=0$ and $a_i=b_i=\tilde{a}_i=\tilde{b}_i=1 \quad
\forall i$, we can obtain it.

\section{Illustrative Examples}
In this section we are going to give illustrative examples\\
\\
\textbf{Example 3.1.} Solve the nearly pentadiagonal linear system
of size $10$ given by
\begin{eqnarray}\left[\begin{array}{cccccccccc}
{3} & {-1} & {3} & {5} & {0} & {0} & {0} & {0} & {0} & {0} \\
{-2} & {2} & {1} & {2} & {0} & {0} & {0} & {0} & {0} & {0} \\
{3} & {-4} & {5} & {5} & {1} & {0} & {0} & {0} & {0} & {0} \\
{0} & {3} & {-2} & {1} & {1} & {3} & {0} & {0} & {0} & {0} \\
{0} & {0} & {6} & {1} & {2} & {5} & {1} & {0} & {0} & {0} \\
{0} & {0} & {0} & {3} & {-3} & {2} & {7} & {-5} & {0} & {0} \\
{0} & {0} & {0} & {0} & {-8} & {1} & {12} & {3} & {-4} & {0} \\
{0} & {0} & {0} & {0} & {0} & {2} & {5} & {3} & {1} & {20} \\
{0} & {0} & {0} & {0} & {0} & {0} & {3} & {11} & {21} & {3} \\
{0} & {0} & {0} & {0} & {0} & {0} & {-2} & {4} & {-9} & {31}
\end{array}\right] \left[\begin{array}{c} {x_{1} } \\
{x_{2} } \\ {x_{3} } \\ {x_{4} } \\ {x_{5} } \\ {x_{6} } \\ {x_{7} }
\\ {x_{8} } \\ {x_{9} } \\ {x_{10} } \end{array}
\right]=\left[\begin{array}{l} {30} \\ {13} \\ {35} \\ {27} \\ {69}
\\ {18} \\ {38} \\ {280} \\ {328} \\ {247}
\end{array}\right]
\end{eqnarray}
by using the \textbf{KNPENTA} algorithm and \textbf{KSNPENTA} algorithm.\\
\\
\textbf{Solution}\\

(i) The application of the \textbf{KNPENTA} algorithm gives:\\
\begin{itemize}
    \item $c_1=3 $(Step 1).\\
    \item $a_n=5, b_1=-2, e_1=-1, f_2=\frac{-2}{3},$ and $c_2=\frac{4}{3} $(Step 3).\\
    \item $e_2=3, e_n=\frac{16}{3}$ and $\tilde{a}_2=\frac{16}{3} $ (Step 5).\\
    \item $[f_3,f_4,f_5,f_6,f_7,f_8,f_9]=[\frac{-9}{4}, -1, \frac{-253}{35}, \frac{-105}{184}, \frac{4012}{2271}, \frac{-368}{101}, \frac{204567}{186433}]$,\\
         $[c_3,c_4,c_5,c_6,c_7,c_8,c_9]=[\frac{35}{4}, 1, \frac{552}{35}, \frac{757}{92}, \frac{-1313}{1514}, \frac{14341}{303}, \frac{4112262}{186433}]$ and \\
          $[e_3,e_4,e_5,e_6,e_7,e_8,e_9]=[12, 2,\frac{934}{35}, \frac{1393}{184}, \frac{26873}{2271}, \frac{-1371}{101}, \frac{-3532041}{186433}] $(Step 6).\\
    \item $f_1=\frac{-91736}{186433}, f_n=\frac{-1203361}{4112262}$ and $c_n=\frac{701215}{19866}$ (Step 7).\\
    \item $[z_1, z_2]=[30,33]$ (Step 8).\\
    \item $[z_3,z_4,z_5,z_6,z_7,z_8,z_9]=[\frac{317}{4}, 32, \frac{8609}{35}, \frac{11475}{184}, \frac{238883}{4542}, \frac{138311}{303}, \frac{129996}{14341}]$ (Step 9).\\
    \item $z_10=\frac{3506075}{9933}$ (Step 10).\\
    \item $[x_1,x_2,x_3,x_4,x_5,x_6,x_7,x_8,x_9,x_{10}]=[1,2,3,4,5,6,7,8,9,10]$ (Step 11).\\
    \end{itemize}
Also the determinant of the matrix $A$ is
$det(A)=-145151505$ by using (2.5).\\
\\
(ii) The application of the \textbf{KSNPENTA} algorithm gives:\\
\hspace*{2.3cm}$X:=nearly\_penta(\tilde{b},b,d,a,\tilde{a},y)=[1, 2, 3, 4, 5, 6, 7, 8, 9, 10]$.\\
\\
\textbf{Example 3.2.} Solve the nearly pentadiagonal linear system
of size $10$ given by
\begin{eqnarray}\left[\begin{array}{cccccccccc}
{0} & {-1} & {3} & {5} & {0} & {0} & {0} & {0} & {0} & {0} \\
{-2} & {2} & {1} & {2} & {0} & {0} & {0} & {0} & {0} & {0} \\
{3} & {-4} & {5} & {5} & {1} & {0} & {0} & {0} & {0} & {0} \\
{0} & {3} & {-2} & {1} & {1} & {3} & {0} & {0} & {0} & {0} \\
{0} & {0} & {6} & {1} & {2} & {5} & {1} & {0} & {0} & {0} \\
{0} & {0} & {0} & {3} & {-3} & {2} & {7} & {-5} & {0} & {0} \\
{0} & {0} & {0} & {0} & {-8} & {1} & {12} & {3} & {-4} & {0} \\
{0} & {0} & {0} & {0} & {0} & {2} & {5} & {3} & {1} & {20} \\
{0} & {0} & {0} & {0} & {0} & {0} & {3} & {11} & {21} & {3} \\
{0} & {0} & {0} & {0} & {0} & {0} & {-2} & {4} & {-9} & {31}
\end{array}\right] \left[\begin{array}{c} {x_{1} } \\
{x_{2} } \\ {x_{3} } \\ {x_{4} } \\ {x_{5} } \\ {x_{6} } \\ {x_{7} }
\\ {x_{8} } \\ {x_{9} } \\ {x_{10} } \end{array}
\right]=\left[\begin{array}{l} {27} \\ {13} \\ {35} \\ {27} \\ {69}
\\ {18} \\ {38} \\ {280} \\ {328} \\ {247}
\end{array}\right]
\end{eqnarray}
by using the \textbf{KNPENTA} algorithm and \textbf{KSNPENTA} algorithm.\\
\\
\textbf{Solution}\\

(i) The application of the \textbf{KNPENTA} algorithm gives:\\
\hspace*{2.3cm}The method is broken down since $c_1=d_1=0$.\\
\\
\hspace*{0.5cm}(ii) The application of the \textbf{KSBPENTA} algorithm gives:\\
\hspace*{2.3cm}$X:=nearly\_penta(\tilde{b},b,d,a,\tilde{a},y)$=$[\frac{-4092987}
{4589918x-4092987},3\,{\frac {1490963
\,x-2728658}{4589918\,x-4092987}},{\frac {11371544\,x-12278961}{
4589918\,x-4092987}},\\
\hspace*{7.6cm}6\,{\frac {3279301\,x-2728658}{4589918\,x-4092987
}},3\,{\frac {7767312\,x-6821645}{4589918\,x-4092987}},1/3\,{\frac {
90274465\,x-73673766}{4589918\,x-4092987}},\\
\hspace*{7.6cm}1/3\,{\frac {95213875\,x-
85952727}{4589918\,x-4092987}},1/5\,{\frac {188851771\,x-163719480}{
4589918\,x-4092987}},\\
\hspace*{7.6cm}1/15\,{\frac {612846296\,x-552553245}{4589918\,x-
4092987}},2/15\,{\frac
{342051296\,x-306974025}{4589918\,x-4092987}}]
_{x=0}$\\
\hspace*{7.5cm}=[1, 2, 3, 4, 5, 6, 7, 8, 9, 10].\\
Also the determinant of the matrix $A$ is $det(A)=61394805$ and for
more details about how to call this procedure, see appendix 1.

\section{Conclusion}
 The methods described here are very effective, provided that optimal
 LU factorization is used. Our symbolic algorithm is competitive with the
 other methods for solving a nearly pentadiagonal linear system
 which appears in many applications.\\
\section{Acknowledgement}

I like to thank Prof. Dr. M. E. A. El-Mikkawy for his valuable
comments and suggestions.\\
\\
\\
\textbf{Appendix 1. } A Maple procedure for solving a nearly pentadiagonal linear system\\
\\
\noindent $> \#$ A Maple Procedure.\\
$> \#$ Written by Dr. A. A. Karawia 18-6-2008.\\
$> \#$ To compute the solution of A general nearly pentadiagonal
Linear system.\\
$>$ restart:\\
$>$ with(linalg,vector,vectdim):\\
$>$ nearly$\_$penta:=proc(bb::vector,b::vector,
d::vector,a::vector,aa::vector,y::vector)\\
 local i,j,k,n; global T,e,c,f,z,X;\\
 n:=vectdim(d):e:=array(1..n):c:=array(1..n):
 f:=array(1..n):z:=array(1..n):X:=array(1..n):\\
$\#$components of the vectors e, c, and f $\#$\\
 c[1]:=d[1]:if c[1]=0 then c[1]:=x; d[1]:=x;fi:
 e[1]:=a[1]:f[2]:=simplify(b[2]/c[1]):\\
 c[2]:=simplify(d[2]-e[1]*f[2]):if c[2]=0 then c[2]:=x;fi:\\
 e[2]:=simplify(a[2]-f[2]*aa[1]):e[n]:=simplify(aa[2]-f[2]*a[n]):
 aa[2]:=e[n]:\\
 \hspace*{1cm}for i from 3 to n-1 do\\
           \hspace*{2cm}f[i]:=simplify((b[i]-bb[i]*e[i-2]/c[i-2])/c[i-1]):\\
           \hspace*{2cm}if i=3 then e[i]:=simplify(a[i]-f[i]*aa[i-1]-bb[i]*a[n]/c[i-2]);
            else e[i]:=simplify(a[i]-f[i]*aa[i-1]);\\
           \hspace*{2cm}fi:\\
            \hspace*{2cm}c[i]:=simplify(d[i]-bb[i]*aa[i-2]/c[i-2]-e[i-1]*f[i]);
           if c[i]=0 then c[i]:=x; fi:\\
 \hspace*{1cm}end do:\\
 f[1]:=simplify((bb[n]-b[1]*e[n-3]/c[n-3])/c[n-2]):\\
 f[n]:=simplify((b[n]-b[1]*aa[n-3]/c[n-3]-f[1]*e[n-2])/c[n-1]):\\
 c[n]:=simplify(d[n]-f[1]*aa[n-2]-e[n-1]*f[n]):if c[n]=0 then c[n]:=x;
 fi:\\
 $\#$ To compute the vector Z $\#$\\
 z[1]:=y[1]:z[2]:=y[2]-f[2]*z[1]:i:='i':\\
 \hspace*{1cm}for i from 3 to n-1 do\\
\hspace*{2cm}     z[i]:=simplify(y[i]-bb[i]*z[i-2]/c[i-2]-f[i]*z[i-1]):\\
 \hspace*{1cm}end do:\\
 z[n]:=simplify(y[n]-f[n]*z[n-1]-f[1]*z[n-2]-b[1]*z[n-3]/c[n-3]):\\
 $\#$ To compute the Solution of the system X. $\#$\\
 X[n]:=z[n]/c[n]:i:='i':\\
 X[n-1]:=simplify((z[n-1]-e[n-1]*X[n])/c[n-1]):\\
 \hspace*{1cm}for i from n-2 by -1 to 2 do\\
\hspace*{2cm}     X[i]:=simplify((z[i]-e[i]*X[i+1]-aa[i]*X[i+2])/c[i]):\\
 \hspace*{1cm}end do:\\
 X[1]:=simplify((z[1]-e[1]*X[2]-aa[1]*X[3]-a[n]*X[4])/c[i]):\\
 $\#$ To compute the determinant T $\#$\\
 T:=subs(x=0,simplify(product(c[r],r=1..n)));\\
 eval(X):\\
end:\\
$> \#$ Call no. 1 for the procedure nearly$\_$penta. $\#$\\
$>$ x:='x':\\
$>$ aa:=aa:=vector([3,2,1,3,1,-5,-4,20]);\\
\hspace*{3cm} aa := [3, 2, 1, 3, 1, -5, -4, 20]\\
$>$a:=vector([-1,1,5,1,5,7,3,1,3,5]);\\
\hspace*{3cm}a:=[-1,1,5,1,5,7,3,1,3,5]\\
$>$ d:=vector([3,2,5,1,2,2,12,3,21,31]);\\
\hspace*{3cm} d:=[3,2,5,1,2,2,12,3,21,31]\\
$>$ b:=vector([-2,-2,-4,-2,1,-3,1,5,11,-9]);\\
\hspace*{3cm} b:=[-2,-2,-4,-2,1,-3,1,5,11,-9]\\
$>$ bb:=vector([0,0,3,3,6,3,-8,2,3,4]);\\
\hspace*{3cm} bb:=[0,0,3,3,6,3,-8,2,3,4]\\
$>$ y:=vector([30,13,35,27,69,18,38,280,328,247]);\\
\hspace*{3cm} y:=[30,13,35,27,69,18,38,280,328,247]\\
$>$ X:=nearly$\_$penta(bb,b,d,a,aa,y);\\
\hspace*{3cm} X := [1, 2, 3, 4, 5, 6, 7, 8, 9, 10]\\
$>$ T;\\
\hspace*{3cm}          -145151505\\
$> \#$ End of call no. 1. $\#$\\
\\
$> \#$ Call no. 2 for the procedure nearly$\_$penta. $\#$\\
$>$ x:='x':\\
$>$ aa:=vector([3,2,1,3,1,-5,-4,20]);\\
\hspace*{3cm} aa := [3, 2, 1, 3, 1, -5, -4, 20]\\
$>$a:=vector([-1,1,5,1,5,7,3,1,3,5]);\\
\hspace*{3cm}a:=[-1,1,5,1,5,7,3,1,3,5]\\
$>$ d:=vector([0,2,5,1,2,2,12,3,21,31]);\\
\hspace*{3cm} d:=[0,2,5,1,2,2,12,3,21,31]\\
$>$ b:=vector([-2,-2,-4,-2,1,-3,1,5,11,-9]);\\
\hspace*{3cm} b:=[-2,-2,-4,-2,1,-3,1,5,11,-9]\\
$>$ bb:=vector([0,0,3,3,6,3,-8,2,3,4]);\\
\hspace*{3cm} bb:=[0,0,3,3,6,3,-8,2,3,4]\\
$>$ y:=vector([27,13,35,27,69,18,38,280,328,247]);\\
\hspace*{3cm} y:=[27,13,35,27,69,18,38,280,328,247]\\
$>$ X:=nearly$\_$penta(bb,b,d,a,aa,y);\\
\hspace*{3cm} X := $[\frac{-4092987} {4589918x-4092987},3\,{\frac
{1490963 \,x-2728658}{4589918\,x-4092987}},{\frac
{11371544\,x-12278961}{ 4589918\,x-4092987}}, 6\,{\frac
{3279301\,x-2728658}{4589918\,x-4092987 }},3\,{\frac
{7767312\,x-6821645}{4589918\,x-4092987}},\\\hspace*{3.2cm}1/3\,{\frac
{ 90274465\,x-73673766}{4589918\,x-4092987}}, 1/3\,{\frac
{95213875\,x- 85952727}{4589918\,x-4092987}},1/5\,{\frac
{188851771\,x-163719480}{ 4589918\,x-4092987}}, 1/15\,{\frac
{612846296\,x-552553245}{4589918\,x-
4092987}},\\\hspace*{3.2cm}2/15\,{\frac
{342051296\,x-306974025}{4589918\,x-4092987}}]$\\
$>$ T;\\
\hspace*{3cm}          61394805\\
$>$ x:=0:X:=map(eval,op(X));\\
\hspace*{3cm} X := [1, 2, 3, 4, 5, 6, 7, 8, 9, 10] \\
$> \#$ End of call no. 2. $\#$

\end{document}